\documentclass{amsart}

\usepackage{amsfonts,amssymb,amsmath,amsthm}
\usepackage[dvips, final]{graphics}
\usepackage{epsfig}
\usepackage{url}

\newtheorem{Definition}{Definition}
\newtheorem{Lemma}{Lemma}
\newtheorem{Remark}{Remark}
\newtheorem{Theorem}{Theorem}
\newtheorem{Proposition}{Proposition}
\newtheorem{Corollary}{Corollary}
\newtheorem{Example}{Example}

\theoremstyle{definition}
\newtheorem*{Ack}{Acknowledgment}

\newcommand{\R}{\mathbb R}

\newcommand{\N}{\mathbb N}

\newcommand{\id}{\operatorname{id}}

\newcommand{\phase}[2]{\operatorname{G_{#1}({#2})}}

\title[Generating Function of Analytical Poisson Structures]{The Universal Generating Function of Analytical Poisson Structures}

\author[B.~Dherin]{Benoit Dherin}
\address{D-MATH, ETH-Zentrum, CH-8092 Z\"urich, Switzerland}
\email{dherin@math.ethz.ch}
\thanks{B. D. acknowledges partial support of SNF Grant No.~21-65213.01}

\begin{document}

   \maketitle

\begin{abstract}
The notion of generating functions of  Poisson structures was first studied in \cite{CDF2005}.
They are special functions which induce, on open subsets of $\R^d$, a Poisson 
structure together with the local symplectic groupoid integrating it. 
A universal generating function was provided in terms of a formal
power series coming from Kontsevich star product. The present
article proves that this universal generating function 
converges for analytical Poisson structures and compares
the induced local symplectic groupoid with the phase space
of Karasev--Maslov.
\end{abstract}


\section*{Introduction}\label{ch2:intro}

Let $U$ be an open subset of $\R^d$ and let $S$ be a function defined on a neighborhood  of $\{0\}\times \{0\}\times U$ 
in $(\R^{d})^*\times (\R^{d})^*\times U$.
In \cite{CDF2005} and in \cite{dherin2004}, it was shown that if $S$  
satisfies  the SGS conditions and the SGA equation 
	\footnote{ The SGA equation may be seen as the equation of an 
	associative product in the operad of exact Lagrangian submanifolds 
	of $(T^*\R^d)^k\times T^*\R^d$ where the compositions are given by symplectic reduction.  
	See \cite{CDF2005-2} for a presentation of this operad.} 
(see Section \ref{GenFunc}), then S generates a Poisson structure on $U$ given by
$$\alpha(x) := \frac{\partial^2 S}{\partial p_1\partial p_2} (0,0,x) $$
 together with the local symplectic groupoid integrating it. Its structure
maps, defined on a neighborhood $G_S(U)$ of $U$ in $T^*U$, are given by
\[
\begin{array}{cccc}
\epsilon(p,x) & = & (0,x)  &\textbf{unit map}\\
i(p,x)      & = & (-p,x) &\textbf{inverse map}\\
s(p,x)      & = & \nabla_{p_2}S(p,0,x) & \textbf{source map}\\
t(p,x)      & = & \nabla_{p_1}S(0,p,x) & \textbf{target map}.
\end{array}
\]
Suppose now that we have a Poisson structure $\alpha$ on a manifold $M$ and a local chart $U$ with
coordinates $x_1,\dots, x_d$. We may pose the reciprocal problem of finding
a generating function on $U$ such that 
$\alpha(x) := \frac{\partial^2 S}{\partial p_1\partial p_2} (0,0,x) $.
Thus, solving this problem for $\alpha$ provides an explicit symplectic realization (in the sense
of \cite{weinstein1983}) of the Poisson manifold $(\alpha,U)$  as well as an explicit 
description of the unique (up to isomorphism) local symplectic groupoid (see \cite{weinstein1987}) 
of $(\alpha,U)$.
As a consequence, observe that knowing the generating function $S$ of $\alpha$ on $U$ allows us
to transform any Hamilton--Poisson system on $U$
$$\dot X(t) = \alpha(X(t))\nabla f(X(t)),\quad X(0) = \xi\in U,$$
where $f\in C^\infty(U)$ into the following integrable Hamilton system on $G_S(U)\subset T^*U$
\begin{eqnarray*}
\dot p(t) & = &  -\nabla_x H_f(p(t),q(t))\\
\dot q(t) & = &  \nabla_p H_f(p(t),q(t)),
\end{eqnarray*}
with Hamiltonian $H_f = -f\circ t$.  The $d$ constant of motions are given by the
components of the source map $s^1,\dots,s^d$.
In this case, results concerning integrable systems may be directly transfered to
Hamilton--Poisson systems.

In fact, for any Poisson structure $\alpha$ and any local chart $U$, we may find a universal
formal generating function for $\alpha$. This is the main result
of \cite{CDF2005}. This universal generating function turns out to be 
the semi-classical part of Kontsevich star product associated to $\alpha$ on $U$,
$$S_\epsilon(\alpha)(p_1,p_2,x) = (p_1+p_2)x +\sum_{n=1}^\infty\frac{\epsilon^n}{n!}\sum_{\Gamma\in T_{n,2}}
W_\Gamma \hat B_\Gamma(p_1,p_2,x).$$
This relates then symplectic groupoids with star products\footnote{
Since then, Karabegov in \cite{karabegov2004-2} has shown 
how to associate a formal symplectic groupoid to a general star product on
a Poisson manifold. His approach differs essentially from ours by the fact
that he defines a formal symplectic groupoid in terms of the
functions on it.} as already foreseen by Coste, Dazord, and Weinstein in \cite{CDW1987}. 
This universal formula is however not a ``true'' generating function in the sense that
it is a formal power series in a formal parameter $\epsilon$. The main result of
this article is the following Theorem. It states that the universal generating function 
converges for analytical Poisson structures and compares the induced local symplectic groupoid
with the construction of Maslov--Karasev \cite{KM1993}.

\begin{Theorem}
Consider a Poisson manifold $M$ with Poisson bivector $\alpha$ 
such that for any local chart $U$ $$\big| \partial^\beta \alpha^{ij}(x)\big|\leq M_U(x)^{|\beta|+1},$$
where $\beta= (\beta_1,\dots,\beta_d)\in\N^d$ is a multi-index and $M_U(x)$ a positive function.

Then, the universal generating function 
$$S_\epsilon(\alpha)(p_1,p_2,x) = (p_1+p_2)x +\sum_{n=1}^\infty\frac{\epsilon^n}{n!}\sum_{\Gamma\in T_{n,2}}
W_\Gamma \hat B_\Gamma(p_1,p_2,x),$$
converges absolutely for $\epsilon \in [0,1)$ and $$\|p_i\| \leq \frac1{64eM_U^2(x)},\quad i = 1,2.$$
Moreover, the local symplectic groupoid induced by $S_\epsilon$ on $U$ is exactly the same -- not only isomorphic --
as the one constructed by Karasev and Maslov in \cite{KM1993}.
In particular, if we consider the Poisson-Hamilton system  
$$\dot X(t) = \alpha(X(t)) p.$$
for small $ p\in (\R^d)^*$, then the end-points of a solution $X$  satisfy
$X(0) = s( p,  q)$ and $X(1) = t( p,  q)$ for 
$ q  = \nabla_{p_2}S(- p, p,X(0))  = \nabla_{p_1}S( p, -p,X(1))$.
\end{Theorem}
The proof of this Theorem is split between Propositions  \ref{Convergence}, \ref{source-target} and \ref{Comparison}.
Note that the identification with Karasev local symplectic groupoid on $U$ allows us to patch all the induced
groupoids in the local charts into a global manifold $G(U)$, which is the local groupoid integrating
$(M,\alpha)$. The question of defining a global generating function will be addressed elsewhere.

\subsection*{Organization of the article}
In Section \ref{LieSystem}, we briefly review standard results on Poisson manifolds and local symplectic groupoids.
Section \ref{GenFunc} introduces the notion of generating functions of a Poisson structures and provides examples.
In Section \ref{UnivSol}, we prove the convergence of the universal
generating function for analytical Poisson structures. In Section \ref{comparison}, we compare the local
symplectic groupoid generated by this convergent generating function and
the one constructed by Maslov and Karasev in \cite{KM1993}.

The results presented here are part of the author PhD thesis \cite{dherin2004}
and may be seen as  a natural development of \cite{CDF2005}.

\begin{Ack}
The author thanks Giovanni Felder and Alberto S. Cattaneo for their help and encouraging support as well
as Alan Weinstein for useful suggestions.
\end{Ack}

	\section{Lie system and local symplectic groupoids}\label{LieSystem}

In this section, we present a system of non-linear PDEs on an open
subset $U$ of $\R^d$ (that we called the reduced Lie system) and we briefly describe how a solution
of this system generates a local symplectic groupoid structure over $U$.
We also explain that any solution of the reduced Lie system may be used to
transform a Poisson-Hamilton system on $U$ into an integrable Hamilton system on
the local symplectic groupoid.
The results stated here were established by  Coste, Dazord, Karasev,
Maslov and Weinstein in \cite{weinstein1983}, \cite{CDW1987}, 
\cite{karasev1987}, \cite{karasev1989} and \cite{KM1993}. Therefore, we
refer the reader to these articles for details and proofs.

\subsection*{The reduced Lie system}

Let $U$ be an open subset of $\R^d$ and $V$ be an open neighborhood of $U$
in $T^*U$. We say that the functions $s^i:V\rightarrow \R$, $i=1,\dots,d$
are a solution of {\bf the reduced Lie system} if for all $i,j=1,\dots,d$
\begin{eqnarray}\label{SourceEq}
\sum_{u=1}^d \frac{\partial s^i}{\partial p^u}\frac{\partial s^j}{\partial x_u}
-\frac{\partial s^i}{\partial x^u}\frac{\partial s^j}{\partial p_u}
& = &  \alpha^{ij}(s^1,\dots,s^d),\label{source}
\end{eqnarray}
for some functions $\alpha^{ij}:U\rightarrow\R$, $i,j=1,\dots,d$
and the $s=(s^1,\dots,s^d)$ satisfy the {\bf initial condition} $s(0,x) = x$
and the {\bf non-degeneracy condition} $$\det \big(\nabla_x s(p,x)\big)\neq 0.$$

\begin{Remark}-\label{conseq}

\begin{enumerate}
\item A necessary condition for the existence of a solution is that the matrix 
$\alpha(x) = \big(\alpha^{ij}(x)\big)_{i,j=1}^d$ is a Poisson structure on $U$

\item If we denote by $\{,\}_\alpha$ the Poisson bracket on $U$ associated to
$\big(\alpha^{ij}(x)\big)_{i,j=1}^d$ and by $\{,\}_\omega$ the Poisson bracket 
associated to the canonical symplectic form $\omega = dp\wedge dq$ on $V\subset T^*U$,
then the reduced Lie system may be written as
$$\{s^i,s^j\}_\omega = \alpha^{ij}(s), \quad i,j=1,\dots,d.$$
It implies that $s = (s^1,\dots,s^d)$ is a Poisson map from
$(V,\omega)$ to $(U,\alpha)$, i.e, that 
$$\{s^*f,s^*g\}_\omega = s^*\{f,g\}_\alpha.$$
for all $f,g\in C^\infty(U)$.
Such a map is also called  {\bf a symplectic realization} of the Poisson manifold
$(U,\alpha)$.

\item Suppose that $s$ is a solution of the reduced Lie system. The initial 
and non-degeneracy conditions implies, by the implicit function theorem,
that there exists an inverse $Q$ such that
$$Q(p, s(p,q)) = q \qquad\textrm{ and }\qquad s(p,Q(p,x)) =x.$$
\end{enumerate}
\end{Remark}

\subsection*{Integrable Hamiltonian lift}

Let  $\alpha$ be a Poisson structure on $U\subset\R^d$ and let $s$ be a solution
of the reduced Lie system. An easy consequence of the fact that $s$ is a Poisson
map from $(V,\omega)$ to $(U,\alpha)$ is the following.
The Hamilton-Poisson system on $U$ with Hamiltonian
$f\in C^\infty(U)$, i.e., 
$$\dot X(t) = \alpha^{ij}(X(t))\partial_jf(X(t)),\quad X(0) = \xi\in U,$$
may be lifted to the Hamilton system on $V$ 
\begin{eqnarray*}
\dot p(t) & = &  -\nabla_x H_f(p(t),q(t))\\
\dot q(t) & = &  \nabla_p H_f(p(t),q(t)),
\end{eqnarray*}
with Hamiltonian $H_f = f\circ s$.
In other words, if $(p(t),q(t))$ is the solution of the Hamilton system on $V$,
then $X(t) := s(p(t),q(t))$ is the solution of the Poisson system on $U$
provided that $s(p(0),q(0)) = \xi$. 
This Hamilton system is in fact an integrable one. Let us construct the $d$ integrals
of the motion as done in \cite{karasev1989}.
Consider on $U$ the time-dependent Hamiltonian 
$K(x,t) =  \langle Q(-pt, x),p\rangle.$
Denote by $\phi_K^t$ the flow generated by $K$ and define
$$t(p,q) := \phi_K^t(s(p,q))\big|_{t=1}.$$
In \cite{karasev1989}, it is shown that $t=(t^1,\dots,t^d)$
is the only map such that 
\begin{eqnarray}
\{t^i,t^j\}_\omega(p,x) & = &  -\alpha^{ij}(t(x,p)),\label{target}\\
\{s^i,t^j\}_\omega(p,x) & = &  0\label{commutation}\\
\quad t(0,x)=x  &\textrm{ and }&\det \big(\nabla_x t(p,x)\big)\neq 0.
\end{eqnarray}
Thus, we may easily prove from (\ref{commutation}) that $t=(t^1,\dots,t^d)$
constitute $d$ constants of the motion of the Hamilton system on $V$.
The system of PDEs defined by equations \eqref{SourceEq}, \eqref{target} and \eqref{commutation}
together with the respective initial and non--degeneracy conditions for 
$s$ and $t$ is called the {\bf Lie system}.

\subsection*{Local groupoid structure}
We may now describe the  the local symplectic groupoid over $U$ generated by
a solution of the Lie system.
Suppose that $(s,t)$  is a solution of the Lie system defined on a open neighborhood
$\phase U{s,t}$ of $U$ in $T^*U$. Thus, $\phase U{s,t}$ inherits from $T^*U$ the canonical symplectic form $dp\wedge dq$
turning it into a symplectic manifold. This will be the local symplectic groupoid over $U$
associated to a solution of the Lie system. The structure maps are namely given by 
$s$ for the {\bf source}, by $t$ for the {\bf target} and by $\epsilon (x) = (0,x)$ for the {\bf unit}
$\epsilon:U\rightarrow \phase U{s,t}$.
The crucial point is to define the product and the inverse in the following way. 
Take $(p,x)\in \phase U{s,t}$. If $p$ is small enough, it is always possible
to find a time-dependent Hamiltonian $f_\tau \in C^\infty(U)$, $\tau\in\R_+$, such that the flow $\Psi_H^\tau$ of its Hamiltonian lift 
$H = -f_\tau \circ t \in C^\infty(\phase U{s,t})$ carries $(0,s(p,x))$ to $(p,x)$ in time $\tau=1$. 
Then, the {\bf inverse} of $(p,x)$ is defined as 
\begin{eqnarray}\label{inv}i(p,x) = \big(\Psi_H^\tau\big)_{\big|\tau=1}^{-1}\big(0,t(p,x)\big),\end{eqnarray}
and for any $(\bar p,\bar x)$ such that $s(p,x) = t(\bar p,\bar x)$, the {\bf product} between
the two points is defined as
\begin{eqnarray}
m\big((p,x),(\bar p,\bar x)\big) & := & \Psi_H^\tau\big((\bar p,\bar x))_{\big|\tau = 1}.\label{mult}
\end{eqnarray}
This local symplectic groupoid is essentially unique,
see for instance \cite{CDW1987} and \cite{karasev1989} for details and proofs.

	\section{Generating functions and local symplectic groupoids}\label{GenFunc}

In this section, we consider special functions which induce on an open subset $U$ of $\R^d$
both a Poisson structure and the local groupoid which integrates it. 
As these functions are in fact usual generating
functions of some Lagrangian submanifolds, we will call them generating functions of the local symplectic groupoid
or of the associated Poisson structure.
We will see how such generating functions provide explicit solutions of the Lie system.
Let us first recall what we mean by generating functions in the context of cotangent bundle.

\begin{Definition}[Generating functions of exact Lagrangian submanifolds]
Consider the cotangent bundle $T^*Q$  over a manifold $Q$ endowed with
its canonical symplectic form $\omega = -d\theta$.Then every $1$-form 
$\mu\in\Omega^1(Q)$ can be considered as an embedding of $Q$ in $T^*Q$
($\mu(Q) \subset T^*Q$). This embedding gives a Lagrangian submanifold
when $\mu^* \omega =0$, i.e when $\mu$ is closed (namely $\mu^* \omega = -\mu^* d\theta = -d\mu^* \theta = -d\mu)$).
We can then associate to each closed $1$-form a Lagrangian submanifold of $T^*Q$. 
Moreover when $\mu$ is exact $\mu = dS$, $S\in C^\infty(Q)$ we call $S$ the generating
function of $\mu(Q)$. 
\end{Definition}
We introduce some notations. Let $U\subset\R^d$ be an open set. Then set
\begin{gather*}
B_n(U)  :=  \big({\R^d}^*\big)^n\times U,\quad B_n^0(U)  :=  \{0\}^n\times U\\
B_n  : = B_n(\R^d),\quad B_n^0  : = B_n^0(\R^d).
\end{gather*}
\begin{Definition}[SGA and SGS]
Let $U\subset\R^d$ be an open subset.  We say that a function $S\in C^\infty(B_2)$ satisfies the
Symplectic Groupoid Associativity (SGA) equation on $U$
if for each $(p_1,p_2,p_3,x)\in B_3(U)$ sufficiently close from $B_3^0(U)$ the following implicit equation for $\bar x,\bar p,\tilde x$ and $\tilde p$,
$$\bar x =\nabla_{p_1}S(\bar p ,p_3,x) ,\quad \bar p =\nabla_{x}S(p_1 ,p_2,\bar x), $$
$$\tilde x =\nabla_{p_2}S(p_1,\tilde p,x) ,\quad \tilde p =\nabla_{x}S(p_2 ,p_3,\tilde x),$$
has an unique solution and if the following additional equation holds
$$S(p_1,p_2,\bar x) + S(\bar p,p_3,x)-\bar x\bar p = S(p_2,p_3,\tilde x)+S(p_1,\tilde p,x)-\tilde x\tilde p .$$
Moreover, if $S(p,0,x)=S(0,p,x)=px$ and $S(p,-p,x) =0,$ we say that $S$ satisfies the
Symplectic Groupoid Structure (SGS) conditions.
At last, a function $S\in C^\infty(B_2)$ satisfying both the SGA equation and the SGS conditions will be called
a  symplectic groupoid generating function on $U$ or simply, when no ambiguities arise, a generating function.
\end{Definition}
Let $S$ be a generating function on $U$. We define the phase space of $S$ on $U$ by 
$$\phase US:=\Big\{(p,x)\in T^*U:\nabla_{p_1}S(0,p,x)\in U \quad\textrm{ and }\quad\nabla_{p_2}S(p,0,x) \in U\Big\}.$$
It is easy to see from the SGS conditions that $\phase US$ is an neighborhood of the zero section in $T^*U$. Thus,
$\phase US$ inherits the canonical symplectic structure of $T^*U$, turning $\phase US$ into a symplectic manifold.
We proved in \cite{CDF2005} and \cite{dherin2004} that the matrix
$$\alpha(x) = \Big(2\nabla_{p_k^1}\nabla_{p_l^2}S(0,0,x)\Big)_{k,l=1}^d$$
is a Poisson structure on $U$. This allows us to call $S$ the {\bf generating function of 
the Poisson structure} $\alpha$.
Moreover, the maps defined by
$$s(p,x)       =  \nabla_{p_2}S(p,0,x)\qquad\textrm{ and }\qquad t(p,x)       =  \nabla_{p_1}S(0,p,x)$$
constitute a solution of the Lie system on $(U,\alpha)$. In particular $s$   is a symplectic realization 
of the Poisson manifold $(U,\alpha)$, i.e., a Poisson map from 
the symplectic manifold $(G_S(U), dp\wedge dq)$ into the Poisson manifold $(G_S(U),\alpha)$.
The local inverses of $s$ and $t$ are given by the explicit formulas
$$Q(p,x) = \nabla_{p_2}S(-p,p,x)\qquad\textrm{ and }\qquad \tilde Q(p,x) = \nabla_{p_1}S(p,-p,x).$$
Finally, summarizing the facts above, we have shown that $G_S(U)$ is the local symplectic 
groupoid integrating $\alpha$ with target $t$, source $s$, unit $\epsilon(x) = (0,x)$ and
inverse $i(p,x)=(-p,x)$. The multiplication space $G_S(U)^{(m)}$  is the graph of
$dS$. 

Let us go through a series of examples. 

\begin{Example}\label{ex:mu}
The function $S(p_1,p_2,x) = \langle x,\mu(p_1,p_2)\rangle$ where $\mu:{\R^d}^*\times {\R^d}^*\rightarrow {\R^d}^*$
is a bilinear associative map on ${\R^d}^*$ always satisfies the SGA equation. We get immediately that 
$$\bar x    =  \langle x,\nabla_{p_1}\mu(\tilde p,p_2)\rangle,\quad \bar p    =  \mu(p_1,p_2),\quad 
\tilde x  =   \langle x,\nabla_{p_2}\mu(p_2,\tilde p)\rangle, \quad \tilde p  =  \mu(p_2,p_3).$$ 
An easy computation shows that  that $S$ satisfies the SGA equation iff $$\mu(p_1,\mu(p_2,p_3))=\mu(\mu(p_1,p_2),p_3).$$
The matrix $$\nabla_{p^1}\nabla_{p^2}S(0,0,x) =\left(2x^i\frac{\partial \mu_i}{\partial p^1_k\partial p^2_l}(0,0,x)\right)_{k,l=1}^d$$
is a Poisson structure if $\mu$ satisfies
\begin{gather*}
\mu(p,0) = \mu(0,p) = p \quad\textrm{ and }\quad \mu(p,-p) = \mu(-p,p) = 0.\\
\end{gather*}

\end{Example}
Notice that it is possible to find a function $S\in C^\infty(B_2)$ which satisfies the SGA equation
but not the SGS conditions. 
\begin{Example}
Consider $T^*\R^2$ and identify ${\R^2}^*$ with the complex plane. We define $S$ as above with $\mu(p_1,p_2) = p_1p_2$ induced
by the complex multiplication. Thus, $$S\big((p_x^1,p_y^1),(p_x^2,p_y^2),(x,y)\big) = 
x(p_x^1p_x^2-p_y^1p_y^2)+y(p_x^1p_y^2+p_y^1p_x^2)$$ satisfies the SGA equation
but not the SGS conditions. By the way the matrix  $$\nabla_{p^1}\nabla_{p^2}S(0,0,x) =\left(
\begin{array}{cc}
x & y  \\
y & -x \\
\end{array} \right)
$$
is not a Poisson matrix.
\end{Example}

\begin{Example}\label{ex:triv}
Let now be $\mu:{\R^d}^*\times {\R^d}^*\rightarrow {\R^d}^*$ the vector addition.
Then, this time, $S(p_1,p_2,x)= x(p_1+p_2)$ satisfies not only the SGA equation but also the SGS conditions.
The associate structure of symplectic groupoid is the trivial one.
The structure maps are $s(p,x) = x$, $t(p,x)=x$, $\epsilon(x) = (0,x)$
$i(p,x) = (-p,x)$. The composition is the fiberwise addition. Moreover the induced Poisson structure on the base
is the null Poisson structure.
\end{Example}

\begin{Example}\label{ex:constant}
Consider $S(p_1,p_2,x) = x(p_1+p_2)+\epsilon \alpha^{ij}p^1_ip^2_j$ where $(\alpha^{ij})_{i,j=1}^d$ is a matrix.
We get that
\begin{gather*}
\bar x^l    = x^l+\epsilon \alpha^{lj}p_j^3, \quad \bar p_l    = p_l^1+p_l^2,\quad
\tilde x^l  = x^l+\epsilon \alpha^{il}p_i^1, \quad  \tilde p_l = p_l^2+p_l^3. 
\end{gather*}
One sees easily that $S$ satisfies the SGA equation for any matrix  $(\alpha^{ij})_{i,j=1}^d$. Moreover,
$$\nabla_{p^1}\nabla_{p^2}S(0,0,x) = (\alpha^{ij})_{i,j=1}^d.$$ Imposing the SGS conditions to $S$ is 
equivalent to say that the matrix $(\alpha^{ij})_{i,j=1}^d$ is skew-symmetric, which also imply that 
$\nabla_{p^1}\nabla_{p^2}S(0,0,x)$ is a Poisson structure of constant rank.
The multiplication space can then be described as 
$$ G^m_S=\Big\{\Big((p_1,x+\epsilon\alpha p_2),(p_2,x-\epsilon\alpha p_1),(p_1+p_2,x)\Big),(p_1,p_2,x)\in B_2\Big\}.$$
The symplectic groupoid structure maps are given by
\begin{eqnarray*}
\epsilon(x) & = & (0,x)\\
i(p,x) & = & (-p,x)\\
s(p,x) &=& x-\epsilon\alpha p\\
t(p,x)& =& x+\epsilon\alpha p.
\end{eqnarray*}
\end{Example}

\begin{Example}\label{ex:CBH}
Suppose ${\R^d}^*$ is given the structure of a Lie algebra, i.e, there is a Lie bracket $[,]:{\R^d}^*\times {\R^d}^*\rightarrow
{\R^d}^*$. In components we have $[p^1,p^2]_k = \alpha_k^{ij}p_i^1p_j^2$ where the $\alpha_k^{ij}$'s are the structure
constants of the Lie algebra. The Baker--Campell--Hausdorff formula 
$$CBH(p_1,p_2)= p_1+p_2+\frac{1}{2}[p_1,p_2]
+\frac{1}{12}\Big([p_1,[p_2,p_2]]+[p_2,[p_2,p_1]]\Big)+\dots,$$ 
induces a structure of group law in a neighborhood of $0$ in $\R^d*$.
As in Example \ref{ex:mu}, define 
$$ S_\epsilon(p_1,p_2,x) = \langle x,\frac{1}{\epsilon}CBH(\epsilon p_1,\epsilon p_2)\rangle. $$
$S$ satisfies the SGA equation because of the associativity of CBH. Moreover, $S$ satisfies also 
the SGS conditions. The induced Poisson structure on $\R^d$ is 
$$\nabla_{p^1}\nabla_{p^2}S(0,0,x) = (\alpha_k^{ij}x^k)_{i,j=1}^d.$$
Notice that, in this case, the induced Poisson structure is linear.
\end{Example}

	\section{The Universal Generating function}\label{UnivSol}

A generating function $S$ produces a Poisson
structure on $U$ together with the local symplectic groupoid integrating it. Here, we investigate the
reciprocal problem. 

\vspace{0.5cm}
{\bf \underline{Question:}
Which  Poisson structures $\alpha$ on $U$ possess a generating function?}
\vspace{0.5cm}

In \cite{CDF2005}, this problem was given a formal answer and a universal
formal generating function was provided. It turns out to be the semi-classical
part of Kontsevich star-product on $U$. In this section, we will first recall
briefly this universal generating function. Then, we will prove that this universal
generating function converges for analytical Poisson structure giving thus a true
generating function, as well as an explicit solution of the Lie system.

Kontsevich, in \cite{kontsevich1997}, parametrizes all deformations (up to gauge equivalences) of the 
usual product of functions on an open set $U\subset  \R^d$ in terms of the (formal) Poisson
structures $\alpha$ (up to gauge equivalences) that one may put on $U$. He writes down the explicit parametrization
$$
B_\epsilon(\alpha)(f,g)(x) = fg(x) + \sum_{n\geq 1} \frac{\epsilon^n}{n!}
\sum_{\Gamma\in G_{n,2}}W_\Gamma B_\Gamma(\alpha)(f,g)(x),
$$
where $G_{n,2}$ is a set of some special graphs, the Kontsevich graphs of type $(n,2)$,
the $W_\Gamma$ are some real numbers associated to them and the $B_\Gamma(\alpha)$ are
special bidifferential operators on $C^\infty(U)$.

The Proposition answering the previous question for an analytical Poisson structure is the following.

\begin{Proposition}\label{Convergence}
Consider a Poisson manifold $M$ with Poisson bivector $\alpha$ 
such that for any local chart $U$ $$\big| \partial^\beta \alpha^{ij}(x)\big|\leq M_U(x)^{|\beta|+1},$$
where $\beta= (\beta_1,\dots,\beta_d)\in\N^d$ is a multi-index and $M_U(x)$ a positive function.

Then, the universal generating function 
$$S_\epsilon(\alpha)(p_1,p_2,x) = (p_1+p_2)x +\sum_{n=1}^\infty\frac{\epsilon^n}{n!}\sum_{\Gamma\in T_{n,2}}
W_\Gamma \hat B_\Gamma(p_1,p_2,x),$$
converges absolutely for $\epsilon \in [0,1)$ and $$\|p_i\| \leq \frac1{64eM_U^2(x)},\quad i = 1,2.$$
\end{Proposition}

The reader may find definitions of Kontsevich graphs, weights and operators in
the original paper \cite{kontsevich1997} and a brief introduction to them 
as well as the definition of Kontsevich trees and symbols in \cite{CDF2005}.
Similar arguments  may be found in \cite{ADS2002} for the case of  linear Poisson structures.

\subsection*{Proof of Proposition \ref{Convergence}}
This follows from a sequence of Lemmas.

\begin{Lemma}
Let $\Gamma\in T_{n,2}$, $p_1,p_2\in B_\rho(0)$. Then,
$$\big|\hat B_\Gamma(p_1,p_2,x)\big| \leq M^{2n-1}\rho^{n+1}.$$
\end{Lemma}

\begin{proof}
The estimate follows from that $\hat B_\Gamma$ is a $(n+1)$-homogeneous polynomial in
the $p$ variables and an $n$-homogeneous polynomial in $\alpha$. Moreover, the is exactly
$n-1$ derivatives differentiating the $\alpha$s.
\end{proof}

\begin{Lemma}
$|T_{n,2}|\leq (16e)^nn!$
\end{Lemma}

\begin{proof}
The graphs in $T_{n,2}$ are exactly the graphs in $G_{n,2}$ such that $\Delta(\Gamma)$ is a
tree, i.e., a graph without cycle. Thus, a graph $\Gamma\in T_{n,2}$ may be described by specifying
$\Delta(\Gamma)$, $n^{n-2}$ choices, by giving each edge of $\Delta(\Gamma)$ an orientation, $2^{n-1}$ choices,
by specifying for each vertex of $\Delta(\Gamma)$ if it has two outgoing edges going to $\{\bar 1,\bar 2\}$,
just one or none, $4^n$ choices, and, at last, by giving the labeling $1$ or $2$ for each edge, $2^n$ choices.
This procedure of counting is by far not optimal. We count too many graphs and also graphs which are
even not Kontsevich graphs. However, each graph in $T_{n,2}$ is counted one times. It gives
us then the very crude estimate,
\begin{eqnarray*}
|T_{n,2}| & \leq & n^{n-2}2^{n-1}4^n2^n\\
&\leq & n^n16^n,
\end{eqnarray*}
which finishes the proof by remarking that $n^n\leq n!e^n$.

\end{proof}

\begin{Definition}
Let $\Gamma\in G_{n,2}$. An aerial vertex $v\in V_\Gamma^a$ is called \textbf{a terminal vertex}
if it has no incoming edge. A terminal vertex $v\in V_\Gamma^a$ is called \textbf{terminal vertex of
type $1$} if one of $\gamma^1(v),\gamma^2(v)$ is in $\{\bar 1,\bar 2\}$ and the other is an aerial vertex.
A terminal vertex $v\in V_\Gamma^a$ is called \textbf{terminal vertex of type $2$} if both $\gamma^1(v)$ and
$\gamma^2(v)$ are aerial vertices.
Trivially, there can be no terminal vertex $v$ such that $\gamma^1(v), \gamma^2(v)\in \{\bar 1,\bar2\}$.
\end{Definition}

\begin{Lemma}
Let  be $\Gamma\in T_{n,2}$. Then, $\Gamma$ has at least one terminal vertex.
\end{Lemma}

\begin{proof}
Suppose $\Gamma$ has no terminal vertex. Let us construct
the following sequence of vertices. Take $v_0\in V_\Gamma^a$. As $v_0$ is not a terminal
vertex, we may chose for $v_1$ any vertex such that $(v_1,v_0)\in E_\Gamma$. As there is
no terminal vertex, we can repeat this procedure infinitely many times. The result is an
infinite sequence of vertices $\{v_l\}_{l\geq0}$  of $V_\Gamma^a$. Now, as $\Gamma$ is 
a Kontsevich tree, there must be no $i,j$ such that $v_i = v_j$, this is a contradiction 
with the fact that $V_\Gamma^a$ is a finite set.
\end{proof}

\begin{Lemma}\label{altern}
Let be  $\Gamma\in T_{n,2}$. Then, we have the following alternative:
\begin{enumerate}
\item There is a terminal vertex of type $1$, $v\in V_\Gamma^a$, such that 
$\Gamma_{|V_\Gamma\backslash\{v\}}\in T_{n-1,2}$.
\item There is a terminal vertex of type $2$, $v\in V_\Gamma^a$, such that 
$\Gamma_{|V_\Gamma\backslash\{v\}} = \Gamma_k\Gamma_l$ where $\Gamma_k\in T_{k,2}$ and
$\Gamma_l\in T_{l,2}$, $k+l = n-1$.
\end{enumerate}
\end{Lemma}

\begin{proof}
Take $\Gamma \in T_{n,2}$. Consider $v$ a terminal vertex in $\Gamma$. If $v$ is of type $1$ 
we are finished. If $v$ is not of type $1$, take another terminal vertex. If this other vertex
is of type $1$, we are finished. If there no terminal vertex of type $1$, there must be at
least one terminal vertex $v$ of type $2$. The restriction of $\Delta(\Gamma)$ to $V_\Gamma^a\backslash\{v\}$
separate $\Delta(\Gamma)$ in two disjoint components $\Delta_k(\Gamma)$ and $\Delta_l(\Gamma)$. They are disjoint
otherwise there would be a cycle in $\Delta(\Gamma)$. Now, it is easy to see that $\Delta_k$ is a tree with
$k$ vertices and $\Delta_l$ a tree with l vertices such that $k+l = n-1$(if not $\Delta(\Gamma)$ is not
a tree). Denote by $\Gamma_k$ and $\Gamma_l$ their associated Kontsevich trees. Then one has that
$\Gamma_{|V_\Gamma\backslash\{v\}} = \Gamma_k\Gamma_l$.
\end{proof}

We illustrate the two cases described by Lemma \ref{altern} by a picture.

\begin{figure}[h]
\begin{center}
\includegraphics{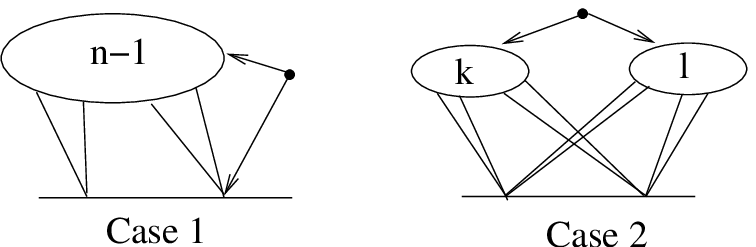}{}
\end{center}
\end{figure}

\begin{Lemma}
Consider  $z\in \mathcal H$ in the hyperbolic upper-half complex plane. Let be $\phi\in[0,2\pi)$ be
an angle. Then, the geometric locus $\mathcal L$ of points $\tilde z\in\mathcal H$ such that $\phi^h(\tilde z, z) = \phi$
is half of the branch of an hyperbola passing though $z$. The other half of the branch
is the locus $ \mathcal  K$ of points such
that $\phi^h(\tilde z,z) = \pi+\phi$.
\end{Lemma}

\begin{proof}
Remark first that for any point $\tilde z$ such that $\phi^h(\tilde z,z) = \phi$ there exist
a point $C_{\tilde z}$ on the real line such that:

(1)The segment joining $z$ to $C_{\tilde z}$  has same length as the segment joining $\tilde z$ to $C_{\tilde z}$.

(2)The angle between the line passing through $C_{\tilde z}$ and $\tilde z$ with the real line is exactly $\phi$.

This is best seen with a drawing:

\begin{figure}[h]
\begin{center}
\includegraphics{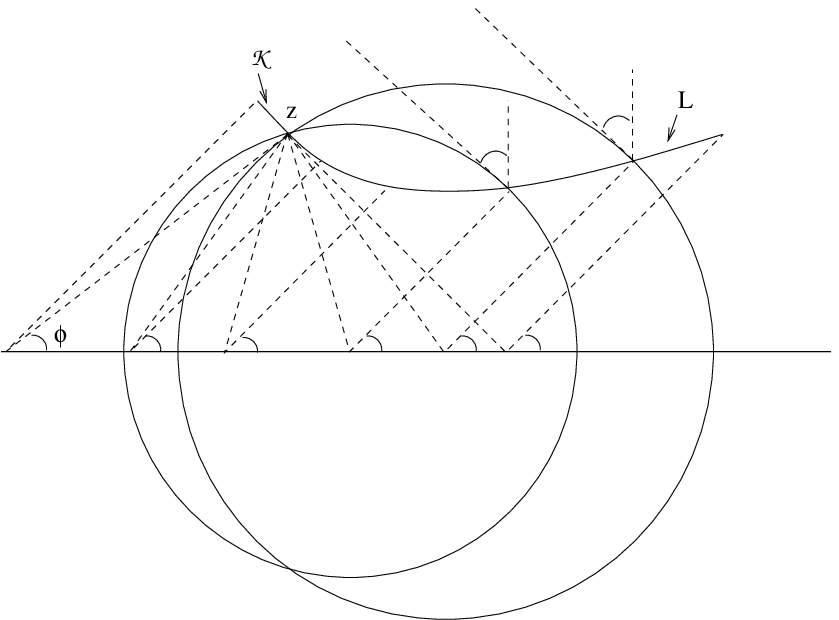}{}
\end{center}
\end{figure}

Translating conditions (1) and (2) into equations we get the equation of a hyperbola passing through $z$. One
sees also easily on the drawing that one half of the branch is the points such that $\phi^h(\tilde z,z) = \phi$
and the other half are the points such that $\phi^h(\tilde z,z) = \pi+\phi$.
\end{proof}

\begin{Corollary}
Let $z\in \mathcal H$ and two angles $\phi_1$ and $\phi_2$. There are at most four couples
$(z_1,z_2)\in\mathcal H^2\backslash D^2 $ such that $\phi^h(z_1,z)= \phi_1$ and $\phi^h(z_2,z) = \phi_2$.
\end{Corollary}

\begin{Lemma}\label{confignum}
For each $\Gamma\in T_{n,2}$ consider the map 
$$\psi_\Gamma:\mathcal H^n\backslash D^n\longrightarrow (S^1)^{2n},$$
defined by
\begin{eqnarray*}
\psi_\Gamma(z_1,\dots,z_n) & = & \Big(\phi^h(z_1,z_{\gamma^1(1)}),\phi^h(z_1,z_{\gamma^2(1)}),\dots\\
&&\dots,\phi^h(z_n,z_{\gamma^1(n)}), \phi^h(z_n,z_{\gamma^2(n)}) \Big).
\end{eqnarray*}
Then, $$|\psi_\Gamma^{-1}(\vec\phi)|\leq 4^n,\quad \vec\phi\in (S^1)^{2n}.$$
\end{Lemma}

\begin{proof}
Given a configuration $\vec z\in \mathcal H^n\backslash D^n$ 
of $\Gamma\in T_{n,2}$ with $$\vec \phi : = \psi_\Gamma(\vec z) = ( \phi_1^1,\phi_1^2,\dots,\phi_n^1,\phi_n^2) ,$$
we have to show that there is at most $4^n$ configurations leading to the same $\vec\phi$.
We do that by induction on $n$. For $n=1$, clearly $|\psi_\Gamma^{-1}(\phi^1,\phi^2)|\leq 4$.
Suppose that it is true for all Kontsevich  trees in $T_{u,2}$ for $u = 1,\dots,n-1$. Then, Lemma \ref{altern},
tells us that we have an alternative. If we are in the first case clearly by induction we have
$$|\phi_\Gamma^{-1}(\vec \phi)| \leq 4^{n-1}.4 = 4^n$$
If we are in the second case, we have that,
$$|\phi_\Gamma^{-1}(\vec \phi)| \leq 4^k4^l4 = 4^n.$$
\end{proof}

\begin{Lemma}
Let $\Gamma\in T_{n,2}$ and $W_\Gamma$ its Kontsevich's weight. Then,
$$|W_\Gamma| \leq 4^n.$$
\end{Lemma}

\begin{proof}
Consider the volume form $$\omega_V = d\phi_1^1\wedge d\phi_1^2\wedge\dots\wedge d\phi_n^1\wedge d\phi_n^2,$$
on $(S^1)^n$. Remark that 
$$W_\Gamma = \frac1{(2\pi)^{2n}}\int_{\mathcal H^n\backslash D^n}\psi_\Gamma^* \omega_V.$$
Thus, by Lemma \ref{confignum}, we have that $$|W_\Gamma|\leq 4^n \frac{V_{(S^1)^{2n}}}{(2\pi)^{2n}},$$
where $V_{(S^1)^{2n}}= (2\pi)^{2n}$ is the volume of $(S^1)^{2n}$. 
\end{proof}

Now, supposing without restricting the generality, that $M_x > 1$,  we can do the following estimate:

\begin{eqnarray*}
\Big|\sum_{n=1}^\infty\epsilon^n\sum_{\Gamma\in T_{n,2}}W_\Gamma \hat B_\Gamma(p_1,p_2,x)\Big| & \leq & 
\sum_{n=1}^\infty \frac{\epsilon^n}{n!}(16e)^n n!4^nM^{2n-1}\rho^{n+1}\\
& \leq & \sum_{n=1}^\infty \epsilon^n \big(64eM^2\big)^n\rho^{n+1}\\
& \leq & \frac1{64eM^2}\sum_{n = 1 }^\infty\epsilon^n,\quad p_1,p_2\in B_\rho(0),
\end{eqnarray*}
which gives the desired result.

	\section{Comparison with Karasev symmetric solution}\label{comparison}

We consider a Poisson manifold $M$ with a Poisson structure $\alpha$ such that in every
chart $U$  for each $x\in U$ there exists a positive function $M_U(x)$ such
that 
 $$\big| \partial^\beta \alpha^{ij}(x)\big|\leq M_U(x)^{|\beta|+1},$$
for all $i,j=1,\dots,d$ and for all multi-index $\beta = (\beta_1,\dots,\beta_d)\in \N^d $.

We may assign  to each chart $U$ the generating function,
$$S_\epsilon(p_1,p_2,x) = x(p_1+p_2)+\sum_{n=1}^\infty \frac{\epsilon^n }{n!}\sum_{\Gamma\in T_{n,2}}W_\Gamma \hat B_\Gamma(p_1,p_2,x)$$
which is convergent, thanks to Proposition \ref{Convergence}, 
for $(p_1,p_2,x)$ taken in a  suitable neighborhood of $B_2^0(U)$.

In particular, there exists a neighborhood $\phase US$ such that $B_1^0(U)\subset \phase U S\subset T^*U$  
and on which $$s(p,q) = \nabla_{p_2}S_\epsilon(p,0,x)\quad\mathrm{and}\quad t(p,q) = \nabla_{p_1}S_\epsilon(0,p,q),$$
take their values in $U$ for all $(p,q)\in \phase US$.

For convenience, set $B(x) := 2\nabla_{p_1}\nabla_{p_2}S(0,0,x) = 2\epsilon \alpha(x)$.
The following Proposition tells us that the points $x_0 = s(p,q)$ and $x_1=t(p,q)$ are
the end points of a curve $x(t)$ in $U$ satisfying a Hamilton-Poisson differential
on $U$ with the linear Hamiltonian $l_p(x) = px$ and initial condition $s(p,q)$.

\begin{Proposition}\label{source-target}
Let be $(\bar p,\bar q)\in \phase US$ . Consider the Poisson
equation $$\dot x^i (t) = B^{ij}(x(t))\bar p_j, \quad x(0) = s(\bar p,\bar q).$$
Then, $$x(1) = t(\bar p,\bar q).$$

Moreover, the Hamilton system associated to $H(p,q) = -t^*l_p(p,q)$, $l_{\bar p}(x) = \langle \bar p,x\rangle,$
on $\phase US$ takes the point $(0,s(\bar p,\bar q))$ to the point $(\bar p,\bar q)$ in time $t = 1$.

\end{Proposition}

\begin{proof}
Consider the Hamilton function $l_{\bar p}(x) = \langle \bar p, x\rangle$ on $U$.  
The vector fields $-\omega^\sharp t^*dl_p$ and $B^\sharp  dl_{\bar p}$ are $t$-related. Thus, the solution 
$(p(t),q(t))$ of the Hamilton system on $\phase US$,
\begin{eqnarray}
\dot p & = & \langle \frac\partial{\partial q}t(p,q),\bar p\rangle\label{pp},\quad p(0) = 0,\\
\dot q & = &-\langle \frac\partial{\partial p}t(p,q),\bar p\rangle,\label{qq}\quad q(0) = s(\bar p,\bar q)
\end{eqnarray}
projects to the solution $x(t) = s(p(t),q(t))$ of the Poisson system
$\dot x = B(x)p$, $x(0) = s(\bar p,\bar q)$.

Remark now that for all $q$ we have that 
\begin{eqnarray*}
\langle \frac{\partial t}{\partial q}(t\bar p,q),\bar p\rangle & = & \Big\langle \nabla_{x}\nabla_{p_1}S_\epsilon(0,t\bar p, q),
\bar p\Big\rangle\\
& = & \bar p+\sum_{i\geq1}\epsilon^i \Big\langle \nabla_{x}\nabla_{p_1}S_\epsilon^{(i)}(t\bar p,0,q),\bar p\Big\rangle\\
& = & \bar p,
\end{eqnarray*}
as for $i\geq 1$, $S_\epsilon^{(i)}$ is a sum over Kontsevich's trees.
Thus, if we insert $p(t) = \bar pt$ in (\ref{pp}) and in (\ref{qq}), Equation
(\ref{pp}) is trivially always satisfied and remains for $q(t)$ the equation
$$\dot q = -\langle\frac{\partial t}{\partial p}(\bar pt, q),\bar p\rangle, \quad
q(0) = s(\bar p,\bar q).$$
For small $\bar p$, there always exists a unique solution for $q(t)$ and as for all $t$ 
we have that $s(p(t),q(t)) = s(\bar p,\bar q)$, then 
$$q(t) = Q(\bar pt,s(\bar p,\bar q)).$$
Thus, $q(1) = Q(\bar p, s(\bar p,\bar q)) = \bar q$ and 
we finally get that $$t(\bar p,\bar q) = t(p(1),q(1)) = x(1)$$.
\end{proof}

In \cite{karasev1987} and \cite{karasev1989}, Karasev gives another local solution of the Lie system on
$U$. Take a point $x\in U$ and a point $ p\in\R^{d*}$ small enough so  that the Poisson system  
$$\dot x^i  =B^{ij}\big(x\big) p_j,\quad x(0) = x$$ has a solution $x(t)$ 
defined for $t\in[-1,1]$. Consider the function
$$Q'(p,x)  = \int_{0}^1 x(t)dt.$$
If $p$ is small enough, we may inverse this function. We define the source of the point $(p,q)$ as
$$s'(p,q)= x,\quad  q = Q'(p,x).$$
The target is defined as $t'(p,x) = s'(-p,x)$.

It happens that $(s',t')$ is a solution of the Lie system in a neighborhood of the null section.
This solution is called \textbf{the symmetric solution}.

The next natural question to ask is how this symmetric solution compares to the one given by the 
generating function.

For that purpose, we will use a Proposition proved in \cite{karasev1987} and 
\cite{karasev1989} by Karasev. Before stating the Proposition, let us make a definition.

\begin{Definition}
Let $(s,t)$ be a solution of the Lie system.
Consider the following differential equation for $p(t)$,
\begin{eqnarray}
\dot p & = & \langle \frac{\partial t}{\partial q}(p,Q(p,x_0)),\bar p\rangle,\quad p(0) = 0,\quad x_0 \in U. \label{exp}
\end{eqnarray}
Then, define the exponential map associated to $(s,t)$ as, $$\exp_{x_0}(\bar p) := p(1).$$
\end{Definition}

\begin{Proposition}[Karasev \cite{karasev1989}]\label{explem}
Let $(s,t)$ and $(s',t')$ be two solutions of the Lie system on $U$ and let $\exp$
and $\exp'$ be their associated exponential maps. Then, the transformation
\begin{eqnarray*}
p' & = & \exp_x'(\exp_x^{-1})(p)\\
q' & = & Q'(p',x),\quad x=s(p,q),
\end{eqnarray*}
is a symplectic map and satisfies
$$s'(p',q') = s(p,q),\quad t'(p',q') = t(p,q).$$
\end{Proposition}

With the help of this Proposition, we may compare the symmetric solution and the one given
by the generating function. 

In fact, in \cite{karasev1989}, Karasev proves that for the symmetric solution $(s',t')$ 
of the Lie system, $\exp_x'(p) = p$. As for the one given by the generating
function, we get from the proof of Proposition \ref{source-target} also that $\exp_x(p) = p$.

\begin{Proposition}[Comparison]\label{Comparison}
The the solution $(s,t)$ of the Lie system given by the symplectic groupoid generating function is exactly the
symmetric solution, i.e., if we consider the Poisson system 
$$\dot x^i (t) = B^{ij}(x)p_j ,\quad x(0) = x,$$
then, $$Q(p,x) = \int_0^1x(t)dt = \nabla_{p_2}S(-p,p,x),$$
$$Q(p,s(p,q)) = q,\quad s(p,Q(p,x)) = x.$$
\end{Proposition}

\begin{proof}
Consider the symmetric solution $(s',t')$. In this case, by Proposition \ref{explem}, the Karasev 
transformation relating $(s,t)$ and $(s',t')$ takes the form
\begin{eqnarray*}
p' & = & p,\\
q' & = & Q'(p,s(p,q)).
\end{eqnarray*}

This transformation is also symplectic. This imposes on the transformation that
$$\frac{\partial q'}{\partial q} = \id\quad\textrm{and}
\quad \frac{\partial q'}{\partial p} = \Big(\frac{\partial q'}{\partial p}\Big)^*.$$
Hence, $$\frac{\partial q'}{\partial q} = \frac{\partial Q'}{\partial x}(p,x)\frac{\partial s}
{\partial q}(p,q) = \id.$$
As we also have that 
$$ \frac{\partial Q}{\partial x}(p,x)\frac{\partial s}{\partial q}(p,q) = \id,$$
then, $$
\frac{\partial Q'}{\partial x}(p,x)
=\frac{\partial Q}{\partial x}(p,x),
$$
for all $p$, $x$. Thus,
\begin{eqnarray}\label{trans}
Q'(p,x) & =  & Q(p,x) + C(p),
\end{eqnarray}
and the Karasev transformation becomes 
\begin{eqnarray*}
p' & = & p,\\
q' & = & Q(p,s(p,q)) + C(p) = q+C(p).
\end{eqnarray*}
Now, 
$$\frac{\partial q'}{\partial p} = \Big(\frac{\partial q'}{\partial p}\Big)^*,$$
implies that 
$$\frac{\partial C}{\partial p}(p) = 
\Big(\frac{\partial C}{\partial p}(p)\Big)^*.$$
Hence, by Stokes, there exists a function $f:\R^{d*}\rightarrow \R$ such that $C(p) = \nabla f(p)$.
Now, as we have for all $p,x$ and $\lambda\in \R$ that,
$$\langle Q(\lambda p,x),p\rangle = \langle x,p\rangle,\qquad
\langle Q'(\lambda p,x),p\rangle = \langle x,p\rangle.$$
equation (\ref{trans}) gives that 
$$\langle C(\lambda p) ,p\rangle = \langle \nabla f(\lambda p),p\rangle = 0.$$
As $$\frac d{d\lambda} f(\lambda p) = \langle \nabla f(\lambda p),p\rangle = 0,$$
we get that $f(\lambda p) = f(0)$, i.e., $f$ is constant and finally
$$\nabla f(p) = c(p) = 0$$ which completes the proof.
\end{proof}

\end{document}